\documentclass{amsproc}
\usepackage{amsmath,amsthm,amsfonts,amscd, amssymb, epsf}
\usepackage{a4wide}
\makeatletter
 \theoremstyle{plain}
\newtheorem{thm}{Theorem}[section]
  \theoremstyle{plain}
  \newtheorem{lem}[thm]{Lemma}
  \newtheorem*{thm*}{Theorem}
\numberwithin{equation}{section}
\usepackage{palatino}
\usepackage{amsthm}
\usepackage{amsfonts}
\usepackage{amscd}
\usepackage{epsf}
\makeatletter

\newcommand{\field}[1]{\ensuremath{\mathbb{#1}}}
\newcommand{\gs}{\Gamma_{1}}
\newcommand{\gb}{\Gamma}
\newcommand{\fs}{\mathfrak{F}}
\newcommand{\fb}{\mathcal{F}}
\newcommand{\rs}{\chi}
\newcommand{\rb}{\psi}

\newcommand{\CC}{\field{C}}
\newcommand{\df}{\equiv}

\newcommand{\HHHH}{\field{H}}

\DeclareMathOperator*{\lto}{\mathnormal{o}(1)}

\DeclareMathOperator{\HH}{\HHHH^3}

\newcommand{\PP}{\field{P}}

\newcommand{\RR}{\field{R}}

\newcommand{\ZZ}{\field{Z}}

\newcommand{\F}{\mathcal{F}}

\newcommand{\scz}{\mathcal{S}}

\newcommand{\D}{\mathcal{D}}
\newcommand{\hil}{\mathcal{H}}

\newcommand{\smat}{\mathfrak{S}}

\DeclareMathOperator{\rep}{Rep}
\DeclareMathOperator{\PSL}{PSL}

\DeclareMathOperator{\R}{Re}
\DeclareMathOperator{\I}{Im}
\DeclareMathOperator{\pc}{PSL(2,\CC)}
\DeclareMathOperator{\GL}{GL}

\DeclareMathOperator{\ldsb}{\mathnormal{ \frac{\phi_{\gb}^{\prime}}{\phi_{\gb}}}}
\DeclareMathOperator{\ldss}{\mathnormal{ \frac{\phi_{\gs}^{\prime}}{\phi_{\gs}}}}
\DeclareMathOperator{\lds}{\mathnormal{ \frac{\phi^{\prime}}{\phi}}}
\DeclareMathOperator{\tr}{tr}

\DeclareMathOperator{\hs}{\mathnormal{\hil(\Gamma,\chi)}}
\DeclareMathOperator{\lp}{\mathnormal{\varDelta}}

\newcommand{\ra}{\rightarrow}

\newcommand{\nsub}{\vartriangleleft}
\makeatother

\begin{document}

\author{Joshua S. Friedman}
\address{ Department of Mathematics and Sciences, United States Merchant Marine Academy, 300 Steamboat Road, Kings Point, NY  11024}
\email{CrownEagle@gmail.com}
\email{friedmanj@usmma.edu}
\email{joshua@math.sunysb.edu}

\title[Analogues of the Artin factorization formula]{Analogues of the Artin factorization formula for the automorphic scattering matrix and Selberg zeta-function associated to a Kleinian group}

\begin{abstract}
For Kleinian groups acting on hyperbolic three-space, we prove factorization formulas for both the Selberg zeta-function and the automorphic scattering matrix. We extend results of Venkov and Zograf from Fuchsian groups, to Kleinian groups, and we give a proof that is simple and extendable to more general groups.
\end{abstract}

\maketitle

\thispagestyle{empty}
\section{Introduction}
In \cite{Venkov3} Venkov and Zograf gave an analogue of Artin's well-known factorization formula. More specifically, they gave factorization formulas for both the Selberg zeta-function and automorphic scattering matrix that are associated to a Fuchsian group (in the context of the Selberg spectral theory of automorphic functions).  They proved the following: 
\begin{thm*}[Venkov-Zograf]
Let $\Gamma$ be a Fuchsian group, $\Gamma_{1}$ a finite-index normal subgroup of $\Gamma,$ then 
\begin{equation} \label{eqIntroFacSZ}
Z(s,\Gamma_{1},\mathbf{1}) = Z(s,\Gamma,U^{\mathbf{1}}) = \prod_{\vartheta \in (\Gamma_{1} \setminus \Gamma)^{*}} Z(s,\gb,\vartheta)^{n_{\vartheta}}. 
\end{equation}
Here $(\Gamma_{1} \setminus \Gamma)^{*}$ is the set of all pairwise inequivalent irreducible unitary representations of the group $\Gamma_{1} \setminus \Gamma;$ for $\vartheta \in (\Gamma_{1} \setminus \Gamma)^{*},$ $n_{\vartheta}$ is the dimension of $\vartheta;$ $Z(s,\Gamma,\vartheta)$ is the Selberg zeta-function; $\mathbf{1}$ is the trivial representation of $\Gamma_{1},$ and $U^{\mathbf{1}}$ is its induced representation to $\Gamma.$ 
\end{thm*} 
In fact, they proved that for any finite-dimensional unitary representation $\chi$ of $\Gamma_{1},$
\begin{equation} \label{introEqSelZetaFac}
Z(s,\Gamma_{1},\chi) = Z(s,\Gamma,U^{\chi}). 
\end{equation}

Venkov and Zograf also proved:
\begin{thm*}[Venkov-Zograf] 
Let $\Gamma$ be a cofinite Fuchsian group, $\Gamma_{1}$ a finite-index normal subgroup of $\Gamma,$ then 
\begin{equation} \label{introEqScatFac}
\phi(s,\Gamma_{1},\chi)\Omega(\Gamma_{1},\chi)^{1-2s} = \phi(s,\Gamma,U^{\chi})\Omega(\Gamma,U^{\chi})^{1-2s}.
\end{equation}  
Here $\phi(s,\,\cdot \, ,\,\cdot \,) $ is the determinant of the automorphic scattering matrix, and  $\Omega(\,\cdot \, ,\,\cdot \,)$ is a constant depending on the group and unitary representation. 
\end{thm*}

In this paper, we give a simple proof of these two theorems---a proof that extends to higher dimensions---for the case of cofinite Kleinian groups acting on hyperbolic three-space. Our proof also applies to the Fuchsian case, and slightly strengthens Equation~\ref{introEqScatFac}: Our proof implies that 
$$\phi(s,\Gamma_{1},\chi) = \phi(s,\Gamma,U^{\chi}) $$ and
$$\Omega(\Gamma_{1},\chi) = \Omega(\Gamma,U^{\chi}).$$

Next we state our results and switch from Fuchsian to Kleinian groups.  
Let $\gb < \pc$ be a cofinite Kleinian group, $\gs$ a finite-index normal subgroup of index $n.$ Let $\rs \in \rep(\gs,V)$ be a finite-dimensional unitary representation of $\gs$ in $V,$ and let $\rb = U^{\rs}\in\rep(\gb,V^{n})$ be its induced representation to $\gb.$

Our main results are as follows:
\begin{thm*} Let $\gb$ be an arbitrary Kleinian group, $\gs$ a finite-index normal subgroup of $\gb.$ Let $\rs \in \rep(\gs,V),$ and $\rb = U^{\rs}.$  Then for $\R{s} > 1,$
\[
Z(s,\gs,\rs)=Z(s,\gb,\rb).\]
\end{thm*}
Here $Z(s,\,\cdot \,,\,\cdot \,)$ is the Selberg zeta-functions (defined in \S\ref{secSZ}).

\begin{thm*}  
Let $\gb$ be an arbitrary Kleinian group, $\gs$ a finite-index normal subgroup of $\gb.$ Then for \mbox{$\R{s} > 1,$}
$$Z(s,\gs,\mathbf{1}) = Z(s,\gb,U^{\mathbf{1}}) = \prod_{\vartheta \in (\gs \setminus \gb)^{*}} Z(s,\gb,\vartheta)^{n_{\vartheta}}. $$
Here $(\gs \setminus \gb)^{*}$ is the set of all pairwise inequivalent irreducible unitary representations of the group $\gs \setminus \gb;$ for $\vartheta \in (\gs \setminus \gb)^{*},$ $n_{\vartheta}$ is the dimension of $\vartheta;$ $Z(s,\gb,\vartheta)$ is the Selberg zeta-function; $\mathbf{1}$ is the trivial representation of $\gs,$ and $U^{\mathbf{1}}$ is its induced representation to $\gb.$
\end{thm*}

Let $\smat(s,\gs,\rs), \smat(s,\gb,\rb)$ be automorphic scattering matrices (defined in \S\ref{secScat}), and set
$$\phi(s,\gs,\rs)  \df \det \smat(s,\gs,\rs),$$
$$\phi(s,\gb,\rb)  \df \det \smat(s,\gb,\rb).$$

\begin{thm*}
Let $\gb$ be a cofinite Kleinian group, $\gs$ a finite-index normal subgroup of $\gb.$ Let $\rs \in \rep(\gs,V),$ and $\rb = U^{\rs}.$ Then for all regular $s\in \CC,$
\[
\phi(s,\gs,\rs)=\phi(s,\gb,\rb).\]
\end{thm*}

I would like to thank Professor Leon Takhtajan for originally suggesting this problem to me, and for reading through this paper. I would also like to thank Peter Zograf for reading this paper and for useful discussions. I would also like to thank the anonymous referee for correcting some errors.

\section{Preliminaries}
In this section we state the preliminary results that will be needed. Our main references are \cite{Friedman2,Friedman1}, \cite{Elstrodt}, and \cite{Venkov}, and \cite{Hejhal1,Hejhal2}. Unless stated otherwise, throughout this section $\Gamma < \pc $ is a cofinite Kleinian group and $\chi \in \rep(\Gamma,V)$ is a finite-dimensional unitary representation of $\Gamma$ in $V.$

\subsection{Cofinite Kleinian Groups}
Let $\Gamma < \pc $ be a cofinite Kleinian group acting on hyperbolic three-space $\HH.$ Let $V$ be a finite-dimensional complex inner product space with inner-product $\left<~,~\right>_{V},$ and let $\rep(\Gamma,V)$ denote the set of finite-dimensional unitary representations of $\Gamma$ in $V.$ Let $\F \subset \Gamma$ be a fundamental domain for the action of $\Gamma$ in $\HH.$ 

Let $\chi \in \rep(\Gamma,V).$ The Hilbert space of $\chi-$\emph{automorphic} functions is the set of measurable functions    
\begin{multline*}
\hs \df  \{ f: \HH \ra V ~|~ f(\gamma P) = \chi(\gamma) f(P)~\forall \gamma \in \Gamma,  P \in \HH, \\ $ and $ \left<f,f \right> \df \int_{\F} \left<f(P),f(P)\right>_V\,dv(P) < \infty \}. 
\end{multline*}
Finally, let   $ \lp = \lp(\Gamma,\chi) $ be the corresponding positive self-adjoint Laplace-Beltrami operator on $\hs.$ 

Next we briefly define the concept of a \emph{singular} unitary representation. Let $\PP$ be the the boundary of $\HH, $ the Riemann sphere. For every $\zeta \in \PP^1$  let  $\Gamma_\zeta $ denote the stabilizer subgroup of  $\zeta$ in $\Gamma, $ 
$$
\Gamma_{\zeta} \df \{~ \gamma \in \Gamma ~ | ~ \gamma \zeta = \zeta ~\},
$$ 
and  let $ \Gamma_\zeta^\prime $ be the maximal torsion-free parabolic subgroup of $\Gamma_\zeta $ (the maximal parabolic subgroup of $\Gamma_\zeta $ that does not contain elliptic elements). A point $\zeta \in \PP^1$ is called a \emph{cusp} of $\Gamma $ if  $\Gamma_\zeta^\prime$ is a free abelian group of rank two.  Two cusps $\zeta_1, \zeta_2$ are $\Gamma-$equivalent  if $\zeta_1 \in \Gamma \zeta_2,$ that is their $\Gamma-$orbits coincide.    

Every cofinite Kleinian group has finitely many equivalence classes of cusps, so we fix a set 
$ \{ \zeta_\alpha \}_{\alpha = 1}^{\kappa(\Gamma)} $ of representatives of these equivalence classes.  For notational convenience we set \mbox{$\Gamma_\alpha \df \Gamma_{\zeta_\alpha} $} and  $\Gamma_\alpha^\prime \df \Gamma_{\zeta_\alpha}^\prime. $  

For each cusp $\zeta_\alpha$ fix an element $B_\alpha \in \pc, $  a lattice $\Lambda_\alpha = \ZZ \oplus \ZZ \tau_\alpha,~\I(\tau_\alpha) > 0, $ and a root of unity $\epsilon_\alpha $ of order 1,2,3,4,or 6 with the following conditions being satisfied:

(1) $\zeta_{\alpha} = B_{\alpha}^{-1} \infty,$

(2) $$  B_{\alpha}\Gamma_{\alpha}^\prime B_{\alpha}^{-1} =  
\left\{ \,\left.\left(\begin{array}{cc}
1 &  b\\
0 &   1\end{array}\right)\,\right|\, b \in \Lambda_\alpha \,\right\},  $$

(3) $$  B_{\alpha}\Gamma_{\alpha} B_{\alpha}^{-1}= \left\{ \,\left.\left(\begin{array}{cc}
\epsilon & \epsilon b\\
0 & \epsilon^{-1}\end{array}\right)\,\right|\, b \in \Lambda_{\alpha},~\epsilon~\text{is some power of $\epsilon_\alpha$} \right\} /\{\pm I\}$$
The group $  B_{\alpha}\Gamma_{\alpha}^\prime B_{\alpha}^{-1} $ acts on  $\CC $ via the lattice $\Lambda_\alpha.$ See \cite[Theorem 2.1.8]{Elstrodt} for more details.

For each cusp  $\zeta_\alpha$ of $\Gamma, $ define the \emph{singular} space 
$$
V_{\alpha} \df \{ v \in V \, | \, \chi(\gamma)v=v, \, \, \, \forall \gamma \in \Gamma_{\alpha}\,\}, 
$$  
and the \emph{almost singular} space 
$$  V_{\alpha}^\prime \df \{ v \in V \, | \, \chi(\gamma)v=v, \, \, \, \forall \gamma \in \Gamma_{\alpha}^\prime \,\}, 
$$ 
where $1 \leq \alpha \leq \kappa(\Gamma). $

A representation   $\chi \in \rep $  is  \emph{singular}  at the cusp $ \zeta_\alpha $ of $\Gamma $ iff the subspace 
$V_{\alpha}  \neq \{0 \}.$ If a cusp is not singular, it is called \emph{regular.}     

For each cusp $\zeta_\alpha,$ set $k_\alpha = \dim_\CC V_\alpha, $   and
$$k(\Gamma,\chi) \df \sum_{\alpha = 1}^{\kappa(\Gamma)} k_\alpha.  $$ 

\subsection{Selberg Zeta-Function} \label{secSZ}
In this section we define the Selberg zeta-function $Z(s,\Gamma,\chi)$ associated to $\Gamma$ and $\chi \in \rep(\Gamma,\chi).$ We allow $\Gamma$ to be an \emph{arbitrary} Kleinian group. See \cite{Friedman2} and \cite[Sections  5.2,5.4]{Elstrodt}, for more details.

Suppose $T\in\Gamma$ is loxodromic (we consider hyperbolic elements
as loxodromic elements). Then $T$ is conjugate in $\pc$ to a unique
element of the form \[
D(T)=\left(\begin{array}{cc}
a(T) & 0\\
0 & a(T)^{-1}\end{array}\right)\]
 such that $a(T)\in\CC$ has $|a(T)|>1$. Let $N(T)$ denote the \emph{norm}
of $T,$ defined by \[
N(T)\df|a(T)|^{2},\]
 and let by $\mathcal{C}(T)$ denote the centralizer of $T$ in $\Gamma.$
There exists a (primitive) loxodromic element $T_{0},$ and a finite
cyclic elliptic subgroup $\mathcal{E}$ of order $m(T),$ generated
by an element $E_{T}$ such that \[
\mathcal{C}(T)=\langle T_{0}\rangle\times\mathcal{E}.\]
 Here $\langle T_{0}\rangle=\{\, T_{0}^{n}~|~n\in\ZZ~\}.$ Next, Let
$\mathfrak{t}_{1},\dots,\mathfrak{t}_{n},$ and $\mathtt{t'_{1}},\dots,\mathtt{t'_{n}}$
denote the eigenvalues of $\chi(T_{0})$ and $\chi(E_{T})$ respectively.
The elliptic element $E_{T}$ is conjugate in $\PSL(2,\CC)$ to an
element of the form \[
\left(\begin{array}{cc}
\zeta(T_{0}) & 0\\
0 & \zeta(T_{0})^{-1}\end{array}\right),\]
 where here $\zeta(T_{0})$ is a primitive $2m(T)$-th root of unity.

For $\R(s)>1$ the Selberg zeta-function $Z(s,\Gamma,\chi)$ is defined
by \[
Z(s,\Gamma,\chi)\df\prod_{\{ T_{0}\}\in\mathcal{R}}~\prod_{j=1}^{\dim V}\prod_{\substack{l,k\geq0\\
c(T,j,l,k)=1}
}\left(1-\mathfrak{t}_{j}a(T_{0})^{-2k}\overline{a(T_{0})^{-2l}}N(T_{0})^{-s-1}\right).\]
 Here the product with respect to $T_{0}$ extends over a maximal
reduced system $\mathcal{R}$ of $\Gamma$-conjugacy classes of primitive
loxodromic elements of $\Gamma.$ The system $\mathcal{R}$ is called
reduced if no two of its elements have representatives with the same
centralizer. See \cite{Elstrodt} Section 5.4 for more details. The
function $c(T,j,l,k)$ is defined by \[
c(T,j,l,k)\df\mathtt{t'_{j}}\zeta(T_{0})^{2l}\zeta(T_{0})^{-2k}.\]
 In \cite{Friedman2}, we gave the meromorphic continuation of $Z(s,\Gamma,\chi)$
to the left half plane under certain technical assumptions.%
\footnote{We assumed that $\Gamma$ was cofinite and had only one class of cusps. We also showed
that the presence of \emph{cuspidal} elliptic fixed points can cause
the zeta-function to have branch points.%
}

A practical way of understanding the zeta-function is via its logarithmic
derivative. For \mbox{$\R(s)>1,$} \begin{equation}
\frac{d}{ds}\log Z(s,\Gamma,\chi)=\sum_{\{ T\}_{\text{lox}}}\frac{\tr(\chi(T))\log N(T_{0})}{m(T)|a(T)-a(T)^{-1}|^{2}}N(T)^{-s}.\label{eq:logZeta}\end{equation}
 Let $W(s,\Gamma,\chi)=\frac{d}{ds}\log Z(s,\Gamma,\chi)$. Then,
for $\R s>1,$ \begin{equation}
Z(s,\Gamma,\chi)=e^{\int W(s,\Gamma,\chi)\, ds+C}\label{eq:zetaInteg}\end{equation}
 where $C$ is chosen so \[
\lim_{s\ra\infty}Z(s,\Gamma,\chi)=1\]
 will be satisfied.

\subsection{Selberg Theory of $\lp$} \label{secSelThe}
In this section we assume that $\Gamma$ is cofinite, and we state some needed results concerning the Selberg trace formula associated to $\Gamma$ and $\chi \in \rep(\Gamma,V).$ We will not need the full trace formula found in \cite{Friedman2,Friedman1}. Rather, only parts of its proof will be needed. More details can be found in \cite{Friedman2,Friedman1}, \cite{Elstrodt}, and \cite{Venkov}.

For $P = z+rj,~P'=z'+r'j \in \HH$ set 
$$ \delta(P,P') \df \cosh(d(P,P')) = \frac{|z-z'|^{2}+r^{2}+ r'^{2}}{2rr'}, $$ where $d$ is the hyperbolic distance in $\HH.$ 

For $k \in \scz([1,\infty)),$ a Schwartz-class function\footnote{Note that by standard approximation techniques, we can (and will) look at functions that are not in $\scz([1,\infty)).$ What we really require is that \eqref{eq:partPonSum} converges absolutely, and uniformly on compacts subsets of $\HH \times \HH;$ and that $h$ (see \eqref{eqSHC}) decays sufficiently fast so that it can be used in the Selberg trace formula (see \cite{Elstrodt}). In fact, we can start with $h,$   a holomorphic function on  $ \{ s \in \CC \, | \, |\I(s)| < 2+ \delta \}$ for some $\delta > 0,$ satisfying $ h(1+z^2) = O( 1+|z|^2)^{3/2 - \epsilon}) $ as  $|z| \ra \infty;$ and recover the function $k.$  },  set 
$$K(P,P')=k(\delta(P,P')).$$ 
Note that for any $\gamma \in \pc;$ 
\begin{equation} \label{eqPointPair}
K(\gamma P,\gamma P') = K(P,P') \quad  \text{ and }  \quad  K( P,\gamma P') = K(\gamma^{-1}P,P').
\end{equation}

For $\Theta\subset \Gamma,$ $\chi \in \rep(\Gamma,V),$ define
\begin{equation} \label{eq:partPonSum}
K(P,P',\Theta,\psi) \df \sum_{\gamma \in \Theta} \chi(\gamma) K(P,\gamma P').
\end{equation}
The series above converges absolutely and uniformly on compact subsets of
$\HH \times \HH.$

For $\lambda \in \CC,$ $\lambda = 1-s^{2},$ the Selberg--Harish-Chandra transform of $k$ is the function $h,$ defined by 
\begin{equation} \label{eqSHC}
h(\lambda) = h(1-s^2) \df \frac{\pi}{s}
\int_{1}^{\infty}k\left(\frac{1}{2}\left(t+\frac{1}{t}\right)\right)
(t^{s}-t^{-s})\left(t-\frac{1}{t}\right)\,\frac{dt}{t},~~\lambda = 1-s^2.
\end{equation}
In addition, let 
$$ g(x) = \frac{1}{2\pi} \int_{\RR} h(1+t^2)e^{-itx}\,dt. $$

We can start with the function $h$ and work backwards to find $k$ \cite[Chapter 3]{Elstrodt}. The pair $h,g$ is said to be \emph{admissible} if   
$h $ is a holomorphic function on  $ \{ s \in \CC \, | \, |\I(s)| < 2+ \delta \}$ for some $\delta > 0,$ satisfying $ h(1+z^2) = O( 1+|z|^2)^{3/2 - \epsilon}) $ as  $|z| \ra \infty.$

For $v,w \in V $ let $ v \otimes \overline{w}$ be the linear operator in $V,$ defined by,    
$ v \otimes \overline{w}(x) = <x,w>v,  $ where $x \in V.$

\begin{lem}[\cite{Friedman1}]
Let  $k \in \scz([1,\infty)) $ and  $h:\CC\ra\CC$ be the Selberg--Harish-Chandra  Transform
of $k.$ Then 
\begin{multline}
\label{E:kernel expansion}
K(P,Q,\Gamma,\chi)   = \sum_{m \in \D} h(\lambda_m)e_m(P) \otimes 
\overline{e_m(Q)}   \\ +
\frac{1}{4\pi}\sum_{\alpha = 1}^{\kappa(\Gamma)} \sum_{l = 1}^{k_\alpha}
\frac{\left[ \Gamma_\alpha:\Gamma_\alpha^\prime \right]}
{\left| \Lambda_\alpha \right|}
\int_{\RR} h \left( 1+t^2 \right) E_{\alpha l}(P,it) \otimes \overline{E_{\alpha l}(Q,it)} \, dt.
\end{multline}

The sum and integrals converge on compact subsets 
of $\HH \times \HH$. Here $\D$ is an indexing set of the eigenfunctions $e_m$ of $\lp$ with corresponding eigenvalues $\lambda_m,$   $E_{\alpha l}(P,s)$ are the Eisenstein series associated to the  singular cusps of $\Gamma$,   $k_\alpha = \dim_\CC V_\alpha $, and $\left| \Lambda_\alpha \right| $ is the Euclidean area of a fundamental domain for the lattice  $ \Lambda_\alpha. $ If a cusp is regular it is omitted from the sum in \eqref{E:kernel expansion}. 
\end{lem}

Next, we split up $K_\Gamma $ as a sum of two kernels. The first kernel 
$$
H_\Gamma(P,Q) =  \frac{1}{4\pi}\sum_{\alpha = 1}^\kappa  
\sum_{l = 1}^{k_\alpha}
\frac{\left[ \Gamma_\alpha:\Gamma_\alpha^\prime \right]}
{\left| \Lambda_\alpha \right|}
\int_{\RR} h \left( 1+t^2 \right) E_{\alpha l}(P,it) 
\otimes \overline{E_{\alpha l}(Q,it)} \, dt,
$$ is not of Hilbert-Schmidt class, while the second kernel $$
L_\Gamma(P,Q) =  \sum_{m \in \D} h(\lambda_m)e_m(P) \otimes 
\overline{e_m(Q)}
$$  
is of trace class.

Suppose $Y>0$ is sufficiently large. Then for all $A > Y,$ there exists a compact set $\F_{A}\subset\HH$
such that 
\begin{equation} \label{eqCuspSec}
\F \df \F_{A}\cup\F_{1}(A)\cup\cdots\cup\F_{\kappa}(A)
\end{equation}
is a fundamental domain for $\Gamma.$  The sets $\F_{\alpha}(A) $ are cusp sectors (see \cite{Elstrodt} Proposition 2.3.9).  It follows that 
 \begin{multline}   \label{lemSingCancel}
\lim_{A \ra \infty} \left( \int_{\F_A} \tr_V (K_{\Gamma}(P,P))  \,dv(P)
-  \int_{\F_A} \tr_V (H_{\Gamma}(P,P))  \,dv(P)  \right) \\
=  \int_\F \tr_V (L_{\Gamma}(P,P))\,dv(P) = 
\sum_{m \in \D} h(\lambda_m)  < \infty.
\end{multline} 
The infinite sum is absolutely convergent.  

Let $\smat(s)$ denote the \emph{automorphic scattering matrix} associated to $\Gamma$ and $\chi$ (see \cite{Friedman1}), and let $$\phi(s) = \det \smat(s).$$ Upon applying the vector form of  the Maa\ss-Selberg  relations (see \cite{Roelcke}, \cite{Venkov}, \cite{Friedman1}), we obtain 
\begin{multline} 
\label{L:SpecTrace}
\int_{\F_A} \tr_V (H_{\Gamma}(P,P))  \,dv(P) = \\ 
g(0)k(\Gamma,\chi) \log(A) - \frac{1}{4\pi} \int_{\RR}\lds(it)
h(1+t^2)\,dt + \frac{h(1)\tr \smat(0)}{4} +  
\lto_{A \ra \infty}.
\end{multline}
The integral on the right hand side converges absolutely.

\section{Induced Representations}

Let $\gb$ be an \emph{arbitrary} Kleinian group, and let $\gs\vartriangleleft\gb$
be a finite-index normal subgroup of index $n$. Let $\fb,\fs$
be fundamental domains of $\gb,\gs$ respectively, with $[\alpha_{i}]_{i=1}^{n}$
a complete set of representatives for the right-cosets of $\gs\setminus\gb,$
satisfying 
\begin{equation}
\label{eqTile}
\fs=\bigcup_{i=1}^{n}\alpha_{i}(\fb).
\end{equation}
Let $V$ be a finite-dimensional hermitian vector space, and let
$\rs$ be a finite-dimensional unitary representation of $\gs$ in
$V$. Set \[
\overline{\rs}(\gamma)=\left\{ \begin{array}{cc}
\rs(\gamma), & \gamma\in\gs\\
0 & \gamma\notin\gs.\end{array}\right.
\]

Let $\rb\df U^{\rs}$ be the induced representation of $\rs$ from $\gs$
to $\gb.$ More explicitly,

\[
\rb:\gb \mapsto{\GL(V}^{n})\]
 and for $\gamma\in\gb$ and $v_{i}\in V,$ 
\begin{equation}
\rb\left(\gamma\right)\left(\sum_{i=1}^{n}\oplus v_{i}\right)=\sum_{i=1}^{n}\oplus\sum_{j=1}^{n}\overline{\rs}(\alpha_{i}\gamma\alpha_{j}^{-1})v_{j}.\label{eqIndDef}
\end{equation}

It follows from (\ref{eqIndDef}) that for $\gamma\in\gb,$ \begin{equation}
\tr_{V^{n}}\rb(\gamma)=\sum_{i=1}^{n}\tr_{V}\overline{\rs}(\alpha_{i}\gamma\alpha_{i}^{-1}).\label{eqTraInd}\end{equation}

We will need the following result (see \cite{RaySinger} and \cite[Theorem 2.1]{Venkov3}):

\begin{lem} \label{lemIso}
There exists an isometry between the Hilbert spaces $\hil(\gs,\rs)$ and $\hil(\gb,\rb),$ which takes the operator $\lp(\gs,\rs)$ to $\lp(\gb,\rb).$
\end{lem}

We abuse notation and call a set $\Theta\subset\gb$ \emph{normal}
if for all $\gamma\in\gb,$ we have $\gamma\Theta\gamma^{-1}=\Theta.$ 

\begin{lem}
\label{lem:intEqu}Let $\gs\vartriangleleft\gb$ be a finite-index,
normal subgroup of $\gb$ of index $n,$ $\rs\in\rep(\gs,V),$ \\ \mbox{$\rb\df U^{\rs}\in\rep(\gb,V^{n})$,}
and let $\Theta\subset\gb$ be a normal subset. Suppose that \[
\int_{\fb}\tr_{V^{n}}K(P,P,\Theta,\rb)\, dv(P)\]
 converges absolutely. Then

\[
\int_{\fb}\tr_{V^{n}}K(P,P,\Theta,\rb)\, dv(P)=\int_{\fs}\tr_{V}K(P,P,\Theta\cap\gs,\rs)\, dv(P).\]

\end{lem}
\begin{proof}
By (\ref{eqTraInd}), \begin{eqnarray*}
\int_{\fb}\tr_{V^{n}}K(P,P,\Theta,\rb)\, dv(P) & = & \sum_{i=1}^{n}\int_{\fb}\sum_{\gamma\in\Theta}\tr_{V}\overline{\rs}(\alpha_{i}\gamma\alpha_{i}^{-1})K(P,\gamma P)\, dv(P)\\
 & = & \sum_{i=1}^{n}\int_{\fb}\sum_{\substack{\gamma\in\Theta\\
\alpha_{i}\gamma\alpha_{i}^{-1}\in\gs}
}\tr_{V}\rs(\alpha_{i}\gamma\alpha_{i}^{-1})K(P,\gamma P)\, dv(P)\\
 & = & \sum_{i=1}^{n}\int_{\fb}\sum_{\gamma\in\Theta\cap\gs}\tr_{V}\rs(\alpha_{i}\gamma\alpha_{i}^{-1})K(P,\gamma P)\, dv(P)\quad\text{(by normality).}\end{eqnarray*}
 Now, since $\Theta$ and $\gs$ are normal, $\gamma\in\Theta\cap\gs$
implies that $\alpha_{i}\gamma\alpha_{i}^{-1}\in\Theta\cap\gs.$ And
as $\gamma$ goes through each element of $\Theta\cap\gs,$ so does
$\alpha_{i}\gamma\alpha_{i}^{-1}.$ So, for each $i,$ setting $\beta=\alpha_{i}\gamma\alpha_{i}^{-1}$
we obtain \begin{eqnarray}
\sum_{i=1}^{n}\int_{\fb}\sum_{\gamma\in\Theta\cap\gs}\tr_{V}\rs(\alpha_{i}\gamma\alpha_{i}^{-1})K(P,\gamma P)\, dv(P) & = & \sum_{i=1}^{n}\int_{\fb}\sum_{\beta\in\Theta\cap\gs}\tr_{V}\rs(\beta)K(P,\alpha_{i}^{-1}\beta\alpha_{i}P)\, dv(P)\nonumber \\
 & = & \sum_{i=1}^{n}\int_{\fb}\sum_{\beta\in\Theta\cap\gs}\tr_{V}\rs(\beta)K(\alpha_{i}P,\beta\alpha_{i}P)\, dv(P)\label{eq:tmp8343} \\
 & = & \sum_{i=1}^{n}\int_{\alpha_{i}(\fb)}\sum_{\beta\in\Theta\cap\gs}\tr_{V}\rs(\beta)K(Q,\beta Q)\, dv(Q)\label{eq:tmp8383}\\
 & = & \int_{\fs}\sum_{\beta\in\Theta\cap\gs}\tr_{V}\rs(\beta)K(Q,\beta Q)\, dv(Q)\label{eq:tmp8291}\\
 & = & \int_{\fs}\tr_{V}K(P,P,\Theta\cap\gs,\rs)\, dv(P),\nonumber \end{eqnarray}
 where in (\ref{eq:tmp8383}) we set $Q=\alpha_{i}P;$  in \eqref{eq:tmp8343} we used Equation~\ref{eqPointPair}, and we used the fact that $dv(\alpha_{i}Q)=dv(Q);$ and in \eqref{eq:tmp8291} we \emph{tiled} the fundamental domain $\fb$ according to \eqref{eqTile}.
\end{proof}

Lemma \ref{lem:intEqu} works well when the point-pair-invariant $k$
gives rise to an integral kernal of trace-class. However, a careful
look at the proof, shows that we could replace $\fb$ by a \emph{truncated}
fundamental domain $\fb_{A}.$

For the rest of this section, we assume that $\Gamma$ is cofinite. 
Let $Y>0$ be sufficiently large so that for $A > Y,$ $\fb$ decomposes into $\fb=\fb_{A}\cup\fb^{A},$
where $\fb_{A}$ is compact and and $\fb^{A}$ is a union of cusp
sectors \eqref{eqCuspSec}. Since \[
\fs=\bigcup_{i=1}^{n}\alpha_{i}(\fb)\]
 it follows that $\fs=\fs_{A}\cup\fs^{A}$ with $\fs_{A}=\bigcup_{i=1}^{n}\alpha_{i}(\fb_{A})$
and $\fs^{A}=\bigcup_{i=1}^{n}\alpha_{i}(\fb^{A});$ $\fs^{A}$
is also a union of cusp sectors. 

\begin{lem}
Let $\gs\vartriangleleft\gb$ be a finite-index, normal subgroup of
$\gb$ of index $n,$ $\rs\in\rep(\gs,V),$ \\ $\rb\df U^{\rs}\in\rep(\gb,V^{n})$,
and let $\Theta\subset\gb$ be a normal subset. Suppose that \[
\int_{\fb_{A}}\tr_{V^{n}}K(P,P,\Theta,\rb)\, dv(P)\]
 converges absolutely. Then

\[
\int_{\fb_{A}}\tr_{V^{n}}K(P,P,\Theta,\rb)\, dv(P)=\int_{\fs_{A}}\tr_{V}K(P,P,\Theta\cap\gs,\rs)\, dv(P).\]

\end{lem}

\section{Factorization Formula of the Selberg zeta-function}
In this section we prove the analogue of the Venkov-Zograf factorization formula (\cite{Venkov3}) for \emph{arbitrary} Kleinian groups. 

Throughout this section, $\gb$ is an arbitrary Kleinian group and $\gs \nsub \gb$ is a finite-index normal subgroup of index $n.$ 
\begin{thm}
\label{thm:ZetaEqual} Suppose that  $\gs \nsub \gb,$
$\rs\in\rep(\gs,V),$ and that $\rb=U^{\rs}\in\rep(\gb,V^{n}).$
Then for $\R(s) > 1,$
\[
Z(s,\gs,\rs)=Z(s,\gb,\rb).\]
\end{thm}
\begin{proof}
By (\ref{eq:zetaInteg}), it suffices to show that \begin{equation}
W(s,\gs,\rs)=W(s,\gb,\rb).\label{eq:EqaLogDer}\end{equation}
 Let $\Theta=\gb^{\text{lox}}\subset\gb$ be the set of \emph{all}
loxodromic elements of $\gb.$ Note that $\Theta$ is \emph{normal}
since the conjugate of a loxodromic element is loxodromic. For $\R s>1,\,\delta>1,$
set \[
k_{s}(\delta)\df\frac{1}{4\pi}\frac{\left(\delta+\sqrt{\delta^{2}-1}\right)^{-s}}{\sqrt{\delta^{2}-1}}.\]
 For $P,P'\in\HH,$ set $K_{s}(P,P')=k_{s}(\delta(P,P'))$, and, using
(\ref{eq:partPonSum}), define $K_{s}(P,P',\Theta,\rb)$ and \\ \mbox{$K_{s}(P,P',\Theta\cap\gs,\rs)$.}
By \cite[page 185-198]{Elstrodt} and \cite[Lemma 6.3]{Friedman2},
it follows that \[
W(s,\gs,\rs)=\int_{\fs}{\tr_{V}K}_{s}(P,P,\Theta\cap\gs,\rs) \,dv(P),\]
 and $ $ \[
W(s,\gb,\rb)=\int_{\fb}{\tr_{V^{n}}K}_{s}(P,P,\Theta,\rb)\,dv(P).\]
 The result now follows from Lemma~\ref{lem:intEqu}. 
\end{proof}
Compare our proof with the proof given in \cite[Theorem 3.1]{Venkov3}.

If we let $\rs = \mathbf{1},$ the trivial one-dimensional representation, it follows that $\rb = U^{\mathbf{1}}$ can be decomposed into irreducible sub-representations, explicitly: 
$$\rb = \bigoplus_{\vartheta \in (\gs \setminus \gb)^{*}} n_{\vartheta} \vartheta, $$
where $(\gs \setminus \gb)^{*}$ denotes the set of all pairwise inequivalent irreducible unitary representations of the group\footnote{Of course, each $\chi$ must be extended from $\gs \setminus \gb$ to $\gb.$} $\gs \setminus \gb,$ $n_{\vartheta}$ is the dimension of $\vartheta.$ For more details see \cite{Kir}. 

Next, by \eqref{eq:logZeta} and \eqref{eq:zetaInteg} it follows that for $\vartheta_{1}, \vartheta_{2} \in \rep(\gb,V_{i}) ~ (i=1,2),$ 
$$Z(s,\gb,\vartheta_{1} \oplus \vartheta_{2}) = Z(s,\gb,\vartheta_{1})Z(s,\gb,\vartheta_{2}). $$
We have
\begin{thm} Let $\gs \nsub \gb.$ Then for $\R{s}>1,$  
$$Z(s,\gs,\mathbf{1}) = Z(s,\gb,U^{\mathbf{1}}) = \prod_{\vartheta \in (\gs \setminus \gb)^{*}} Z(s,\gb,\vartheta)^{n_{\vartheta}}. $$
\end{thm}

\section{Factorization Formula for the Determinant of the Scattering Matrix} \label{secScat}
In this section we prove the analogous result to Theorem~\ref{thm:ZetaEqual}, for the automorphic scattering matrix. Throughout this section, $\gb$ is a cofinite Kleinian group.

Let $\gs \nsub \gb,$ $\rs \in \rep(\gs,V),$ and $\rb = U^{\rs}.$ Let $k \in \scz([1,\infty)),$ and let $h$ be given by \eqref{eqSHC}. 
It follows form Lemma~\ref{lemIso} that $\lp(\gs,\rs)$ and $\lp(\gb,\rb)$ have the same eigenvalues. Hence (from \S\ref{secSelThe}) 
\[
\int_{\fb}\tr_{V^{n}}L(P,P,\gb,\rb)\, dv(P)=\int_{\fs}\tr_{V}L(P,P,\gs,\rs)\, dv(P) = \sum_{m \in \mathcal{D}}h(\lambda_{m}).
\]
Now, since 
\[
\int_{\fb_{A}}\tr_{V^{n}}K(P,P,\gb,\rb)\, dv(P)=\int_{\fs_{A}}\tr_{V}K(P,P,\gs,\rs)\, dv(P),
\]
it follows from Equation~\ref{lemSingCancel} that

\[
\int_{\fb_{A}}\tr_{V^{n}}H(P,P,\gb,\rb)\, dv(P) = \int_{\fs_{A}}\tr_{V}H(P,P,\gs,\rs)\, dv(P) + \lto_{A \ra \infty}.
\]

However, from  Equation~\ref{L:SpecTrace}, we obtain
\begin{multline} 
\int_{\F_A} \tr_{V^{n}}H(P,P,\gb,\rb)  \,dv(P) = \\ 
g(0)k(\gb,\rb) \log(A) - \frac{1}{4\pi} \int_{\RR}\ldsb(it)
h(1+t^2)\,dt + \frac{h(1)\tr \smat_{\gb}(0)}{4} +  
\lto_{A \ra \infty},
\end{multline}
and
\begin{multline} 
\int_{\F_A} \tr_{V}H(P,P,\gs,\rs)  \,dv(P) = \\ 
g(0)k(\gs,\rs) \log(A) - \frac{1}{4\pi} \int_{\RR}\ldss(it)
h(1+t^2)\,dt + \frac{h(1)\tr \smat_{\gs}(0)}{4} +  
\lto_{A \ra \infty}.
\end{multline}

Hence, we have
\begin{multline} \label{eqScatTrace}
g(0)k(\gb,\rb) \log(A) - \frac{1}{4\pi} \int_{\RR}\ldsb(it)
h(1+t^2)\,dt + \frac{h(1)\tr \smat_{\gb}(0)}{4} \\ = g(0)k(\gs,\rs) \log(A) - \frac{1}{4\pi} \int_{\RR}\ldss(it)
h(1+t^2)\,dt + \frac{h(1)\tr \smat_{\gs}(0)}{4} +  
\lto_{A \ra \infty}.
\end{multline}

Equation~\ref{eqScatTrace} is true for \emph{all} admissible pairs $h,g.$
Hence we must have  
$$k(\gb,\rb) = k(\gs,\rs),  $$
$$\tr \smat_{\gb}(0) = \tr \smat_{\gs}(0),$$ and 
\begin{equation} \label{eqLogScatEq} 
\ldsb(z) = \ldss(z): 
\end{equation}
We will prove \eqref{eqLogScatEq} shortly.
\begin{thm}
Let  $\gs \nsub \gb,$
$\rs\in\rep(\gs,V),$  $\rb=U^{\rs}\in\rep(\gb,V^{n}).$
Let $$\phi_{\gs}(s) \df \det \smat_{\gs}(s) \df \det \smat(s,\gs,\rs),$$
$$\phi_{\gb}(s) \df \det \smat_{\gb}(s) \df \det \smat(s,\gb,\rb).$$
Then for all regular $s,$
\[
\phi_{\gs}(s)=\phi_{\gb}(s).\]
\end{thm}
\begin{proof}
For $r > 0,z \in \CC,$ let $h(z) =e^{-r z^{2}},$  (and it follows) 
$$g(x) = \frac{e^{-r}}{\sqrt{4\pi r}} e^{-x^{2}/(4r)}. $$ 
By considering the asymptotics of \eqref{eqScatTrace}, it follows that 
$$k(\gb,\rb) = k(\gs,\rs),  $$ and 
$$\tr \smat_{\gb}(0) = \tr \smat_{\gs}(0).$$
Hence 
$$ 
\int_{\RR}\ldsb(it) h(1+t^2)\,dt = \int_{\RR}\ldss(it) h(1+t^2)\,dt,
$$
which implies that
$$ \int_{\RR}\ldsb(it) e^{-r t^2}\,dt = \int_{\RR}\ldss(it) e^{-r t^2}\,dt \quad (r > 0).
$$
Now, by the functional equation for $\smat_{\gb}(s)$ (see \cite{Friedman2})
$$\smat_{\gb}(s)\smat_{\gb}(-s) = I. $$ 
So
$$\phi_{\gb}(s) \phi_{\gb}(-s) = 1.$$ Hence  $\ldsb(it)$ is an even function of $t;$ so is $\ldss(it).$ Thus
$$
\int_{0}^{\infty}\ldsb(it) e^{-r t^2}\,dt = \int_{0}^{\infty}\ldss(it) e^{-r t^2}\,dt \quad (r > 0).
$$
Next, the substitution $u = t^{2}$ allows us to rewrite the above integral as a Laplace transform, and by uniqueness (the Laplace transform is invertible), it follows that 
$\ldsb(it) = \ldss(it),$ and (by analytic continuation) 
$$\ldsb(s) = \ldss(s). $$  
Integrating and exponentiating, gives us a constant $C_{1},$ so that
$$\phi_{\gb} = C_{1} \cdot \phi_{\gs}. $$  However, the functional equation (for $\smat_{\gb}(s)$) implies that 
$$\phi_{\gb}(0) = (-1)^{(k_{\gb} -\tr \smat_{\gb}(0) )/2} = (-1)^{(k_{\gs} -\tr \smat_{\gs}(0) )/2} = \phi_{\gs}(0) \neq 0. $$ Hence $C_{1}=1.$
\end{proof}

\bibliographystyle{amsalpha} \bibliographystyle{amsalpha}
\bibliography{vza}

\end{document}